\documentclass{article}[12pt]
\usepackage{amssymb,amsfonts}

\begin{document}

\begin{center}
{\large\bf On the convergence of generalized power series satisfying an algebraic ODE}
\end{center}

\begin{center}
\bf R.\,R.\,Gontsov, I.\,V.\,Goryuchkina
\end{center}
\bigskip

\begin{abstract} We propose a sufficient condition of the convergence of a generalized power series formally satisfying an algebraic (polynomial)
ordinary differential equation. The proof is based on the majorant method.
\end{abstract}

\section{Introduction}

In the paper we study some properties of generalized power series 
\begin{eqnarray}\label{series}
\varphi=\sum_{k=0}^{\infty}c_kx^{s_k}, \qquad c_k\in{\mathbb C}, \quad s_0\prec s_1\prec\ldots\in{\mathbb C}, \quad
\lim_{k\rightarrow\infty}{\rm Re}\,s_k=+\infty,
\end{eqnarray}
formally satisfying an ordinary differential equation
\begin{eqnarray}\label{ADE}
F(x,y,\delta y,\ldots,\delta^n y)=0
\end{eqnarray}
of order $n$ with respect to the unknown $y$, where $F(x,y_0,y_1,\ldots,y_n)\not\equiv0$ is a polynomial of $n+2$ variables and $\displaystyle \delta=x\frac d{dx}$.
Here $\prec$ is a usual ordering by first difference: $s_k\prec s_{k+1}$ iff ${\rm Re}\,s_k<{\rm Re}\,s_{k+1}$ or ${\rm Re}\,s_k={\rm Re}\,s_{k+1}$,
${\rm Im}\,s_k<{\rm Im}\,s_{k+1}$.

Note that substituting the series (\ref{series}) into the equation (\ref{ADE}) makes sense, as only a finite number of terms in $\varphi$
contribute to any term of the form $cx^s$ in the expansion of $F(x,\varphi,\delta\varphi,\ldots,\delta^n\varphi)$ in powers of $x$.
Indeed, $\delta^j\varphi=\sum_{k=0}^{\infty}c_ks_k^jx^{s_k}$ and an equation $s=s_{k_0}+s_{k_1}+\ldots+s_{k_m}$ has a finite number of
solutions $(s_{k_0},s_{k_1},\ldots,s_{k_m})$, since ${\rm Re}\,s_k\rightarrow+\infty$. Furthermore, for any integer $N$ an inequality
$s_{k_0}+s_{k_1}+\ldots+s_{k_m}\preceq N$ has also a finite number of solutions, so powers of $x$ in the expansion of
$F(x,\varphi,\delta\varphi,\ldots,\delta^n\varphi)$ are well ordered with respect to $\prec$.

Earlier in the paper \cite{GS} generalized power series of the form (\ref{series}), with $s_0<s_1<\ldots\in{\mathbb R}$, were studied. There
was proved (without the assumption $s_k\rightarrow+\infty$) that they form a differential ring, and if the series (\ref{series}) satisfies
the equation (\ref{ADE}), then $\lim_{k\rightarrow\infty}s_k=+\infty$. Furthermore, the exponents $s_k\in\mathbb R$ of the formal solution
(\ref{series}) generate a finite $\mathbb Z$-module. 
Here we prove this fact in the case of the complex exponents $s_k\in\mathbb C$ (Lemma 2).

For the generalized power series (\ref{series}) one may naturally define the {\it order}
$$
{\rm ord}\,\varphi=s_0,
$$
and this is also well defined for any polynomial in $\varphi,\delta\varphi,\ldots,\delta^n\varphi$  with coefficients of the form $\alpha\,x^{\beta},$ $\alpha,\;\beta\in\mathbb{C}$.

The main result of the paper is the following sufficient condition of the convergence of (\ref{series}).
\medskip

{\bf Theorem 1.} {\it Let the generalized power series $(\ref{series})$ formally satisfy the equation $(\ref{ADE})$,
$\,\displaystyle\frac{\partial F}{\partial y_n}(x,\varphi,\delta\varphi,\ldots,\delta^n\varphi)\ne0$ and
\begin{eqnarray}\label{ordineq}
{\rm ord}\,\frac{\partial F}{\partial y_j}(x,\varphi,\delta\varphi,\ldots,\delta^n\varphi)\succeq
{\rm ord}\,\frac{\partial F}{\partial y_n}(x,\varphi,\delta\varphi,\ldots,\delta^n\varphi), \quad j=0,1,\ldots,n.
\end{eqnarray}
Then for any sector $S$ of sufficiently small radius with the vertex at the origin and of the opening less than $2\pi$,
the series $\varphi$ converges uniformly in $S$.}
\medskip

This theorem in a somewhat different form has been formulated in \cite[Th. 3.4]{Br} for the case of the {\it real} exponents $s_k\in\mathbb R$.
And in the case of the {\it integer} exponents it becomes Malgrange's theorem \cite{Mal} on the convergence of a formal solution
$\hat\varphi=\sum_{k=0}^{\infty}c_kx^k\in{\mathbb C}[[x]]$ of the equation (\ref{ADE}). An idea of our proof is based on the
construction of a majorant algebraic equation and was used in \cite{GG} for estimating the radius of convergence of a formal
solution $\hat\varphi\in{\mathbb C}[[x]]$ satisfying the conditions of Malgrange's theorem. It originally comes from
\cite[Ch. 1, \S7]{BG}, similar ideas are seemed to be already appeared in \cite{Ca} for studying some properties of divergent
formal solutions $\hat\varphi\in{\mathbb C}[[x]]$ of (\ref{ADE}).

\section{Auxiliary lemmas}

The first auxiliary lemma is a generalization of the corresponding lemma from \cite{Mal}, where it has been proved in the case
of the {\it integer} exponents $s_k$. Here we use technical tools of that paper.
\medskip

{\bf Lemma 1.} {\it Under the assumptions of Theorem 1, there exists an integer $\mu'\geqslant0$ such that for any integer $\mu\geqslant \mu'$ a transformation
$$
y=\sum_{k=0}^{\mu}c_kx^{s_k}+x^{s_{\mu}}u
$$
reduces the equation $(\ref{ADE})$ to an equation of the form
\begin{eqnarray}\label{auxileq}
L(\delta)u=N(x,u,\delta u,\ldots,\delta^nu),
\end{eqnarray}
where $L$ is a polynomial of degree $n$, and
$N$ is a finite sum of monomials of the form
$$
\alpha\,x^{\beta}u^{q_0}(\delta u)^{q_1}\ldots(\delta^nu)^{q_n}, \quad \alpha,\beta\in{\mathbb C},\;{\rm Re}\,\beta>0,\;q_i\in{\mathbb Z}_+.
$$
}

{\bf Proof.} Making a transformation $y=x^{s_0}v$ and taking into consideration the equality $\delta(x^{s_0}v)=x^{s_0}(\delta+s_0)v$, we have the relations
$$
\delta^iy=x^{s_0}(\delta+s_0)^iv, \qquad i=1,\ldots,n,
$$
and therefore come from (\ref{ADE}) to an equation
\begin{eqnarray}\label{eq}
H(x,v,\delta v,\ldots,\delta^nv)=0,
\end{eqnarray}
where the function $H(x,v_0,v_1,\ldots,v_n)$ is a finite sum of monomials of the form
$$
c\,x^{l+ms_0}v_0^{q_0}v_1^{q_1}\ldots v_n^{q_n}, \qquad c\in{\mathbb C}, \quad l,m,q_i\in{\mathbb Z}_+.
$$
We also may assume that $l+ms_0\succeq 0$ (multiplying, if necessary, the equation (\ref{eq}) by a corresponding $x^K$, $K\in{\mathbb N}$). The obtained equation
(\ref{eq}) has a formal solution
$$
\tilde\varphi=\sum_{k=0}^{\infty}c_kx^{s_k-s_0}, \qquad 0\prec s_1-s_0\prec s_2-s_0\prec\ldots\;.
$$

For any $\mu\in{\mathbb Z}_+$, the formal series $\tilde\varphi$ can be represented in the form
$$
\tilde\varphi=\sum_{k=0}^{\mu}c_kx^{s_k-s_0}+x^{s_{\mu}-s_0}\psi=\varphi_{\mu}+x^{s_{\mu}-s_0}\psi, \quad
{\rm ord}\,\psi=s_{\mu+1}-s_{\mu}\succ0.
$$
Then denoting $\Phi=(\tilde\varphi,\delta\tilde\varphi,\ldots,\delta^n\tilde\varphi)=\Phi_{\mu}+x^{s_{\mu}-s_0}\Psi$ and applying the Taylor formula
to the relation $H(x,\Phi)=0$, we have
\begin{eqnarray}\label{Taylor}
0 & = & H(x,\Phi_{\mu}+x^{s_{\mu}-s_0}\Psi)=H(x,\Phi_{\mu})+x^{s_{\mu}-s_0}\sum_{i=0}^n\frac{\partial H}{\partial v_i}(x,\Phi_{\mu})\psi_i+\nonumber \\
  &  & +\frac12x^{2(s_{\mu}-s_0)}\sum_{i,j=0}^n\frac{\partial^2H}{\partial v_i\partial v_j}(x,\Phi_{\mu})\psi_i\psi_j+\ldots,
\end{eqnarray}
where $\psi_i=(\delta+s_{\mu}-s_0)^i\psi$.

One can easily check that the assumption $\displaystyle\frac{\partial F}{\partial y_n}(x,\varphi,\delta\varphi,\ldots,\delta^n\varphi)\ne0$ of Theorem 1 implies $\displaystyle\frac{\partial H}{\partial v_n}(x,\Phi)\ne0$, and the assumption (\ref{ordineq}) implies
$$
{\rm ord}\,\frac{\partial H}{\partial v_i}(x,\Phi)\succeq {\rm ord}\,\frac{\partial H}{\partial v_n}(x,\Phi), \qquad i=0,1,\ldots,n.
$$
Let $\displaystyle\theta={\rm ord}\,\frac{\partial H}{\partial v_n}(x,\Phi)$. Then each formal series
$\displaystyle\frac{\partial H}{\partial v_i}(x,\Phi)$ is of the form
$$
\frac{\partial H}{\partial v_i}(x,\Phi)=b_ix^{\theta}+c_ix^{\theta_i}+\ldots, \qquad i=0,1,\ldots,n, 
$$
where $b_i\in{\mathbb C}$ $(b_n\neq 0)$, $c_i\in{\mathbb C}^*$ and $\theta_i\succ\theta$. Define a polynomial
\begin{equation}
L(\xi)=b_0+b_1(\xi+s_{\mu}-s_0)+\ldots+b_n(\xi+s_{\mu}-s_0)^n\label{Euler}
\end{equation}
of degree $n$ and choose a number $\mu$ such that the following three conditions hold:
$$
1)\, {\rm Re}(s_{\mu}-s_0)>{\rm Re\,}\theta, \qquad 2)\, {\rm Re}(s_{\mu+1}-s_{\mu})>0 \qquad 3)\, L(\xi)\ne0\quad\forall\xi\succ0.
$$
Now we show that such a number $\mu$ is from the statement of the lemma.

Let us note that
\begin{eqnarray*}
{\rm ord}\left(\frac{\partial H}{\partial v_i}(x,\Phi)-\frac{\partial H}{\partial v_i}(x,\Phi_{\mu})\right) & = &
{\rm ord}\biggl(x^{s_{\mu}-s_0}\sum_{j=0}^n\frac{\partial^2 H}{\partial v_i\partial v_j}(x,\Phi_{\mu})\psi_j+\ldots\biggr)=\\
& = & \theta+\alpha_i,\quad{\rm Re}\,\alpha_i>0, \qquad i=0,1,\ldots,n,
\end{eqnarray*}
since ${\rm Re}(s_{\mu}-s_0)>{\rm Re\,}\theta$, and the real parts of the numbers $\displaystyle {\rm ord}\,\frac{\partial^2 H}{\partial v_i\partial v_j}(x,\Phi_{\mu})$ are nonnegative (the same is true for the orders of the other partial derivatives of the function $H$). Therefore,
$$
\frac{\partial H}{\partial v_i}(x,\Phi_{\mu})=b_ix^{\theta}+\tilde c_ix^{\tilde\theta_i}+\ldots, \quad \tilde\theta_i\succeq\theta+\alpha_i.
$$ 
Now from the relation (\ref{Taylor}) and conditions ${\rm Re}(s_{\mu+1}-s_{\mu})>0$, ${\rm Re}(s_{\mu}-s_0)>{\rm Re\,}\theta$, it follows that
$$
{\rm Re}\bigl({\rm ord}\,H(x,\Phi_{\mu})-(s_{\mu}-s_0+\theta)\bigr)>0.
$$
Hence the relation (\ref{Taylor}) can be divided by $x^{s_{\mu}-s_0+\theta}$, and we obtain the equality of the form
$$
L(\delta)\psi-N(x,\psi,\delta\psi,\ldots,\delta^n\psi)=0,
$$
where the polynomial $L$ is defined by the formula (\ref{Euler}), and $N$ is such as in the statement of the lemma.
Thus, the transformation
$$
v=\sum_{k=0}^{\mu}c_kx^{s_k-s_0}+x^{s_{\mu}-s_0}u
$$
reduces the equation (\ref{eq}) to the equation
$$
L(\delta)u=N(x,u,\delta u,\ldots,\delta^nu),
$$
with a formal solution $u=\psi$. $\quad\Box$
\medskip

{\bf Remark 1.} The condition $L(\xi)\ne0$, for all $\xi\succ0$, is not used in the proof of Lemma 1, but we add it from the beginning as it will be used
in the proof of the next lemma and in the sequel.
\medskip

Let us define an additive semi-group $\Gamma$ generated by a (finite) set of power exponents of the variable $x$ containing in 
$N(x,u,\delta u,\ldots,\delta^nu)$, and let $r_1,\ldots,r_{\nu}$ be generators of this semi-group, that is,
$$
\Gamma=\{m_1r_1+\ldots+m_{\nu}r_{\nu} \mid m_i\in{\mathbb Z}_+,\;\sum_{i=1}^{\nu}m_i>0\}, \qquad {\rm Re\,}r_i>0.
$$
The second auxiliary lemma is a consequence of the first one and describes a structure of the set of power exponents $s_k\in\mathbb C$ of the formal solution
(\ref{series}).
\medskip

{\bf Lemma 2.} {\it All the numbers $s_k-s_{\mu}$, $k\geqslant\mu+1$, belong to an additive semi-group $\Gamma$.}
\medskip

{\bf Proof.} We use the fact that the generalized formal power series
$\psi=\sum_{k=\mu+1}^{\infty}c_kx^{s_k-s_{\mu}}$ satisfies the relation
\begin{eqnarray}\label{auxrel}
L(\delta)\psi=N(x,\psi,\delta\psi,\ldots,\delta^n\psi).
\end{eqnarray}
The first term of the left hand side of (\ref{auxrel}) is $L(s_{\mu+1}-s_{\mu})\,c_{\mu+1}x^{s_{\mu+1}-s_{\mu}}$ (note that $L(s_{\mu+1}-s_{\mu})\ne0$, as
$s_{\mu+1}-s_{\mu}\succ0$), while the first term of the right hand side of (\ref{auxrel}) is a monomial $\alpha\,x^{\beta}$, since other monomials
$\alpha\,x^{\beta}\psi^{q_0}(\delta\psi)^{q_1}\ldots(\delta^n\psi)^{q_n}$, $\sum_{i=0}^nq_i>0$, have higher orders. Therefore,
$$
L(s_{\mu+1}-s_{\mu})\,c_{\mu+1}x^{s_{\mu+1}-s_{\mu}}=\alpha\,x^{\beta},
$$
and $s_{\mu+1}-s_{\mu}=\beta\in\Gamma$.

Further, the second term of the left hand side of (\ref{auxrel}) is $L(s_{\mu+2}-s_{\mu})\,c_{\mu+2}x^{s_{\mu+2}-s_{\mu}}$, while the second term of the right hand side of (\ref{auxrel}) is of the form $\alpha\,x^{\beta+m(s_{\mu+1}-s_{\mu})}$, $m\in{\mathbb Z}_+$, $\beta\in\Gamma$. Therefore,
$$
L(s_{\mu+2}-s_{\mu})\,c_{\mu+2}x^{s_{\mu+2}-s_{\mu}}=\alpha\,x^{\beta+m(s_{\mu+1}-s_{\mu})},
$$
and $s_{\mu+2}-s_{\mu}=\beta+m(s_{\mu+1}-s_{\mu})\in\Gamma$. In the analogous way, one obtains that all the numbers
$s_{k}-s_{\mu}$, $k>\mu$, belong to the semi-group $\Gamma$.  {\hfill $\Box$}
\medskip

We may assume that the generators $r_1,\ldots,r_{\nu}$ of the semi-group $\Gamma$ are linearly independent over $\mathbb Z$. This is provided by the following lemma.
\medskip 

{\bf Lemma 3.} {\it There are complex numbers $\rho_1,\ldots,\rho_{\tau}$ linearly independent over $\mathbb Z$, such that all ${\rm Re}\,\rho_i>0$ and an additive 
semi-group $\Gamma'$ generated by them contains the above semi-group $\Gamma$ generated by $r_1,\ldots,r_{\nu}$.}
\medskip

{\bf Proof.} Let $r_1,\ldots,r_{\tau}$ be a maximal system of linearly independent over $\mathbb Z$ elements from the set $\{r_1,\ldots,r_{\nu}\}$. It is sufficient to prove that if we add any number $b$ with ${\rm Re}\,b>0$ such that $r_1,\ldots,r_{\tau},b$ become linearly dependent over $\mathbb Z$, then the semi-group $G$ generated by $r_1,\ldots,r_{\tau},b$ is contained in some semi-group $\Gamma'$ generated by $\tau$ linearly independent over $\mathbb Z$ complex numbers $\rho_1,\ldots,
\rho_{\tau}$ having positive real parts. We may assume that 
\begin{eqnarray}\label{b}
b=m_1r_1+\ldots+m_{\tau-j}r_{\tau-j}-m_{\tau-j+1}r_{\tau-j+1}-\ldots-m_{\tau}r_{\tau}, \qquad m_i\in{\mathbb N}, \quad 1\leqslant j\leqslant\tau-1,
\end{eqnarray}
and we prove the existence of such a semi-group $\Gamma'$ by the induction with respect to the number $j$ of the signes ''$-$'' before the coefficients in the linear combination (\ref{b}). (We assume that all $m_i\ne0$. In general, if some coefficients $m_i$ are equal to zero, then the corresponding generators $r_i$ are included without changes into the set of generators of a new semi-group $\Gamma'$, i.~e., $\rho_i=r_i$ for such $r_i$. In this case the reader can easily make the corresponding changes in the reasonings below.)

For $j=1$ we have
$$
b=m_1r_1+\ldots+m_{\tau-1}r_{\tau-1}-m_{\tau}r_{\tau}.
$$
There are rational positive numbers $p_1/q_1,\ldots,p_{\tau-1}/q_{\tau-1}$ such that $p_1/q_1+\ldots+p_{\tau-1}/q_{\tau-1}=1$ and
\begin{eqnarray}\label{rational}
\frac{p_i}{q_i}<\frac{m_i{\rm Re}\,r_i}{m_{\tau}{\rm Re}\,r_{\tau}}, \qquad i=1,\ldots,\tau-1.
\end{eqnarray}
Indeed, the intersection of the box
$$
\{(x_1,\ldots,x_{\tau-1})\in{\mathbb R}^{\tau-1}\mid 0<x_i<m_i{\rm Re}\,r_i/m_{\tau}{\rm Re}\,r_{\tau}\}
$$
with the hyperplane
$$
\pi=\{x_1+\ldots+x_{\tau-1}=1\}
$$
is a non-empty open subset of $\pi$, since $\sum_{i=1}^{\tau-1}m_i{\rm Re}\,r_i/m_{\tau}{\rm Re}\,r_{\tau}>1$ in view of the condition ${\rm Re}\,b>0$. Hence, we can choose a point in this intersection that has rational coordinates.

Generators $\rho_1,\ldots,\rho_{\tau}$ of a semi-group $\Gamma'$ that are linearly independent over $\mathbb Z$ now can be defined as follows:
\begin{eqnarray*}
\rho_1 & = & r_1-\frac{p_1}{q_1m_1}m_{\tau}r_{\tau}, \quad\ldots,\quad \rho_{\tau-1}=r_{\tau-1}-\frac{p_{\tau-1}}{q_{\tau-1}m_{\tau-1}}m_{\tau}r_{\tau},\\
\rho_{\tau} & = & \frac1{(q_1m_1)\ldots(q_{\tau-1}m_{\tau-1})}r_{\tau}.
\end{eqnarray*}
Then ${\rm Re}\,\rho_i>0$ according to (\ref{rational}). Furthermore $G\subset\Gamma'$, as $b=m_1\rho_1+\ldots+m_{\tau-1}\rho_{\tau-1}$ and
\begin{eqnarray*}
r_1 & = & \rho_1+\frac{p_1}{q_1m_1}m_{\tau}r_{\tau}=\rho_1+n_1\rho_{\tau}, \quad n_1\in{\mathbb N},\\
\ldots & & \ldots \\
r_{\tau-1} & = & \rho_{\tau-1}+\frac{p_{\tau-1}}{q_{\tau-1}m_{\tau-1}}m_{\tau}r_{\tau}=\rho_{\tau-1}+n_{\tau-1}\rho_{\tau}, \quad n_{\tau-1}\in{\mathbb N},\\
r_{\tau} & = & (q_1m_1)\ldots(q_{\tau-1}m_{\tau-1})\rho_{\tau}.
\end{eqnarray*}

For an arbitrary $j>1$ we write the number $b$ in the form $b=b'-m_{\tau}r_{\tau}$, where
$$
b'=m_1r_1+\ldots+m_{\tau-j}r_{\tau-j}-m_{\tau-j+1}r_{\tau-j+1}-\ldots-m_{\tau-1}r_{\tau-1}
$$ 
has $j-1$ signes ''$-$'' in its representation of the form (\ref{b}). Thus, we may apply an inductive assumption to the number $b'$ and write it as follows:
$$
b'=m_1\rho_1+\ldots+m_{\tau-j}\rho_{\tau-j},
$$ 
where the numbers $\rho_1,\ldots,\rho_{\tau-j},\ldots,\rho_{\tau-1}$ are linearly independent over $\mathbb Z$ and expressed via linear combinations of 
$r_1,\ldots,r_{\tau-1}$ with rational coefficients (conversely, $r_1,\ldots,r_{\tau-1}$ are expressed via linear combinations of $\rho_1,\ldots,\rho_{\tau-1}$ with positive integer coefficients). Hence, $r_{\tau}$ and $\rho_1,\ldots,\rho_{\tau-j},\ldots,\rho_{\tau-1}$ are linearly independent over $\mathbb Z$, and we conclude 
for $b=m_1\rho_1+\ldots+m_{\tau-j}\rho_{\tau-j}-m_{\tau}r_{\tau}$ as in the case $j=1$. {\hfill $\Box$}

\section{Proof of Theorem 1}

For the simplicity of the presentation, we give the proof in the case of two generators of the semi-group $\Gamma$:
$$
\Gamma=\{lr_1+mr_2\mid l,m\in{\mathbb Z}_+,\;l+m>0\}, \qquad {\rm Re}\,r_{1,2}>0.
$$
As will be shown further, in this case we deal with functions of two variables. In the case of an arbitrary number $\nu$ of generators, as can be easily seen, 
the proof is analogous, only functions of higher number of variables are involved.

We should prove the convergence of the generalized formal power series
$$
\psi=\sum_{k=\mu+1}^{\infty}c_kx^{s_k-s_{\mu}},
$$
which satisfies the equality
\begin{eqnarray}\label{auxiliary}
L(\delta)\psi=N(x,\psi,\delta\psi,\ldots,\delta^n\psi)
\end{eqnarray}
obtained in Lemma 1. According to Lemma 2, all the exponents $s_k-s_{\mu}$ belong to the semi-group $\Gamma$:
$$
s_k-s_{\mu}=lr_1+mr_2, \qquad (l,m)\in M\subset{\mathbb Z}_+^2\setminus\{0\},
$$
for some subset $M$ such that the map $k\mapsto(l,m)$ is a bijection from ${\mathbb N}\setminus\{1,\ldots,\mu\}$ to $M$. Then
\begin{eqnarray*}
L(\delta)\psi&=&\sum_{k=\mu+1}^{\infty}L(s_k-s_{\mu})\,c_kx^{s_k-s_{\mu}}=\sum_{(l,m)\in M}L(lr_1+mr_2)\,c_{l,m}x^{lr_1+mr_2}=\\
&=&\sum_{(l,m)\in{\mathbb Z}_+^2\setminus\{0\}}L(lr_1+mr_2)\,c_{l,m}x^{lr_1+mr_2}
\end{eqnarray*}
(in the last series one puts $c_{l,m}=0$, if $(l,m)\not\in M$).

Without lost of generality we may assume that all $|s_k-s_{\mu}|=|lr_1+mr_2|\geqslant1$, since $\lim_{k\rightarrow\infty}(s_k-s_{\mu})=\infty$. (In the opposite case
we make the transformation $u=\sum_{k=\mu+1}^{\nu}c_kx^{s_k-s_{\mu}}+w$, where $\nu$ is such that all $|s_k-s_{\mu}|\geqslant1$ for $k>\nu$. This transformation 
reduces the equation (\ref{auxileq}) to an equation of the same form with respect to the unknown $w$, with the formal solution 
$w=\sum_{k=\nu+1}^{\infty}c_kx^{s_k-s_{\mu}}$.)

The function $N(x,u_0,u_1,\ldots,u_n)$ determining the right hand side of the equality (\ref{auxiliary}) is a finite sum of the form
$$
N(x,u_0,u_1,\ldots,u_n)=\sum_{p,q,Q}\alpha_{p,q,Q}\,x^{pr_1+qr_2}u_0^{q_0}u_1^{q_1}\ldots u_n^{q_n},
$$
where $\alpha_{p,q,Q}\in\mathbb C$, $(p,q)\in{\mathbb Z}_+^2\setminus\{0\}$, $Q=(q_0,q_1,\ldots,q_n)\in{\mathbb Z}_+^{n+1}$. Thus the equality (\ref{auxiliary}) is written as follows:
\begin{eqnarray}\label{auxiliaryx}
\sum_{(l,m)\in{\mathbb Z}_+^2\setminus\{0\}}L(lr_1+mr_2)\,c_{l,m}x^{lr_1+mr_2}=\sum_{p,q,Q}\alpha_{p,q,Q}\,x^{pr_1+qr_2}
\psi^{q_0}(\delta\psi)^{q_1}\ldots(\delta^n\psi)^{q_n},
\end{eqnarray}
where
$$
\delta^j\psi=\sum_{k=\mu+1}^{\infty}c_k(s_k-s_{\mu})^jx^{s_k-s_{\mu}}=\sum_{(l,m)\in{\mathbb Z}_+^2\setminus\{0\}}(lr_1+mr_2)^jc_{l,m}x^{lr_1+mr_2}.
$$
Therefore, we have the coincidence of the following two formal power series of two independent variables $z_1, z_2$:
\begin{eqnarray}\label{auxiliaryz}
\sum_{(l,m)\in{\mathbb Z}_+^2\setminus\{0\}}L(lr_1+mr_2)\,c_{l,m}z_1^lz_2^m=\sum_{p,q,Q}\alpha_{p,q,Q}\,z_1^pz_2^q
\psi_0^{q_0}\psi_1^{q_1}\ldots\psi_n^{q_n},
\end{eqnarray}
where
$$
\psi_j=\sum_{(l,m)\in{\mathbb Z}_+^2\setminus\{0\}}(lr_1+mr_2)^jc_{l,m}z_1^lz_2^m\in{\mathbb C}[[z_1,z_2]], \quad j=0,1,\ldots,n.
$$
Indeed, the coefficient $a_{l,m}$ of a monomial $z_1^lz_2^m$ in the right hand side of (\ref{auxiliaryz}) coincides with the coefficient of the corresponding monomial 
$x^{lr_1+mr_2}$ in the right hand side of the equality (\ref{auxiliaryx}), since for each pair $(l,m)\in{\mathbb Z}_+^2\setminus\{0\}$ there is no another pair $(l',m')$ 
such that $l'r_1+m'r_2=lr_1+mr_2$ (in view of the linear independence of the numbers $r_1$, $r_2$ over $\mathbb Z$). Hence, $a_{l,m}=L(lr_1+mr_2)\,c_{l,m}$.

To prove the convergence of $\psi_0\in{\mathbb C}[[z_1,z_2]]$ in some neighbourhood of the origin, we construct an equation
\begin{eqnarray}\label{majorant}
\sigma W=\sum_{p,q,Q}|\alpha_{p,q,Q}|\,z_1^pz_2^qW^{q_0}W^{q_1}\ldots W^{q_n},
\end{eqnarray}
whose right hand side is obtained from that of the equality (\ref{auxiliaryz}) by the change of the coefficients $\alpha_{p,q,Q}$ to their absolute values and all the $\psi_j$ to the one variable $W$. The number $\sigma$ is defined by the formula
$$
\sigma=\inf\limits_{k>{\mu}}\frac{|L(s_k-s_{\mu})|}{|s_k-s_{\mu}|^n}=\inf\limits_{(l,m)\in M}\frac{|L(lr_1+mr_2)|}{|lr_1+mr_2|^n}
$$
and is a positive real number, since $L(\xi)\ne0$ for all $\xi\succ0$, and $\lim\limits_{k\rightarrow\infty}{|L(s_k-s_{\mu})|/|s_k-s_{\mu}|^n}=|b_n|>0$ (recall that 
$L(\xi)=b_0+b_1(\xi+s_{\mu}-s_0)+\ldots+b_n(\xi+s_{\mu}-s_0)^n$, $b_n\ne0$). The equation (\ref{majorant}) possesses a unique holomorphic near the origin solution
$$
W=\sum\limits_{(l,m)\in{\mathbb Z}_+^2\setminus\{0\}}A_{l,m}\,z_1^lz_2^m
$$
satisfying the condition $W(0,0)=0$. This follows from the theorem on implicit function. One can write the coefficients $A_{l,m}$ in the form
$$
A_{l,m}=C_{l,m}|lr_1+mr_2|^n, \qquad |C_{l,m}|\leqslant|A_{l,m}|.
$$
Further we prove that the convergent near the origin power series
$$
\widetilde W=\sum\limits_{(l,m)\in{\mathbb Z}_+^2\setminus\{0\}}C_{l,m}\,z_1^lz_2^m
$$
is majorant for the formal power series $\psi_0$, that is,
$$
C_{l,m}\in{\mathbb R}_+, \qquad |c_{l,m}|\leqslant C_{l,m}\quad\forall (l,m)\in{\mathbb Z}_+^2\setminus\{0\},
$$
which will imply the convergence of $\psi_0$ in some neighbourhood of the origin.

First we use the equality (\ref{auxiliaryz}) to obtain recursive expressions for the coefficients $c_{l,m}$. Denote by $\phi$ the formal power series from
the right hand side of this equality,
$$
\phi=\sum_{p,q,Q}\alpha_{p,q,Q}\,z_1^pz_2^q\psi_0^{q_0}\psi_1^{q_1}\ldots\psi_n^{q_n}\in{\mathbb C}[[z_1,z_2]],
$$
then (\ref{auxiliaryz}) implies
\begin{eqnarray}\label{clm}
L(lr_1+mr_2)\,c_{l,m}=\Bigl.\frac{\partial_{z_1}^l\partial_{z_2}^m\phi}{l!\,m!}\Bigr|_{z_1=z_2=0},
\end{eqnarray}
where $\partial_{z_1}$ is the partial derivative with respect to $z_1$, and $\partial_{z_2}$ is that with respect to $z_2$. To express
$\partial_{z_1}^l\partial_{z_2}^m\phi(0,0)$, let us apply the formulae for the derivation of a product,
\begin{eqnarray*}
\partial_{z_1}^l(f_0\ldots f_n)&=&\sum_{l_0+\ldots+l_n=l}\frac{l!}{l_0!\,l_1!\ldots l_n!}\,\;(\partial_{z_1}^{l_0}f_0)\;\ldots\;(\partial_{z_1}^{l_n}f_n),
\qquad f_i\in{\mathbb C}[[z_1,z_2]],\\
\partial_{z_1}^l\partial_{z_2}^m(f_0\ldots f_n)&=&\sum_{\begin{array}{c}\scriptstyle l_0+\ldots+l_n=l\\ \scriptstyle m_0+\ldots+m_n=m\end{array}}
l!\,m!\,\;\frac{\partial_{z_1}^{l_0}\partial_{z_2}^{m_0}f_0}{l_0!\,m_0!}\;\ldots\;
\frac{\partial_{z_1}^{l_n}\partial_{z_2}^{m_n}f_n}{l_n!\,m_n!}.
\end{eqnarray*}
Thus, we have
\begin{eqnarray}\label{form1}
\Bigl.\frac1{l!\,m!}\,\partial_{z_1}^l\partial_{z_2}^m(z_1^pz_2^q\psi_0^{q_0}\ldots\psi_n^{q_n})\Bigr|_{{z_1}={z_2}=0}=
\sum_{\begin{array}{c}\scriptstyle l_0+\ldots+l_n=l-p\\ \scriptstyle m_0+\ldots+m_n=m-q\end{array}}
\Bigl.\frac{\partial_{z_1}^{l_0}\partial_{z_2}^{m_0}\psi_0^{q_0}}{l_0!\,m_0!}\,\ldots
\frac{\partial_{z_1}^{l_n}\partial_{z_2}^{m_n}\psi_n^{q_n}}{l_n!\,m_n!}\Bigr|_{z_1={z_2}=0}
\end{eqnarray}
for any $l\geqslant p$, $m\geqslant q$ (and $=0$, if $l<p$ or $m<q$). Every partial derivative in this sum is expressed as follows:
\begin{eqnarray}\label{form2}
\Bigl.\frac{\partial_{z_1}^{l_j}\partial_{z_2}^{m_j}\psi_j^{q_j}}{l_j!\,m_j!}\Bigr|_{{z_1}={z_2}=0}&=
&\sum_{\begin{array}{c}\scriptstyle \lambda_1+\ldots+\lambda_{q_j}=l_j\\ \scriptstyle \mu_1+\ldots+\mu_{q_j}=m_j\end{array}}
\Bigl.\frac{\partial_{z_1}^{\lambda_1}\partial_{z_2}^{\mu_1}\psi_j}{\lambda_1!\,\mu_1!}\ldots
\frac{\partial_{z_1}^{\lambda_{q_j}}\partial_{z_2}^{\mu_{q_j}}\psi_j}{\lambda_{q_j}!\,\mu_{q_j}!}\Bigr|_{{z_1}={z_2}=0}=\\ \nonumber
&=&\sum_{\begin{array}{c}\scriptstyle \lambda_1+\ldots+\lambda_{q_j}=l_j\\ \scriptstyle \mu_1+\ldots+\mu_{q_j}=m_j\end{array}}
c_{\lambda_1,\mu_1}(\lambda_1r_1+\mu_1r_2)^j\ldots c_{\lambda_{q_j},\mu_{q_j}}(\lambda_{q_j}r_1+\mu_{q_j}r_2)^j
\end{eqnarray}
(note that this expression is equal to zero, if $l_j+m_j<q_j$). Combining the formulae (\ref{form1}), (\ref{form2}) we obtain
\begin{eqnarray}\label{form3}
\Bigl.\frac{\partial_{z_1}^l\partial_{z_2}^m\phi}{l!\,m!}\Bigr|_{z_1={z_2}=0}=\sum_{p,q,Q}\alpha_{p,q,Q}
\sum_{\begin{array}{c}\scriptstyle l_0+\ldots+l_n=l-p\\ \scriptstyle m_0+\ldots+m_n=m-q\end{array}}
\Bigl.\frac{\partial_{z_1}^{l_0}\partial_{z_2}^{m_0}\psi_0^{q_0}}{l_0!\,m_0!}\,\;\ldots\;
\frac{\partial_{z_1}^{l_n}\partial_{z_2}^{m_n}\psi_n^{q_n}}{l_n!\,m_n!}\Bigr|_{z_1={z_2}=0},
\end{eqnarray}
where
\begin{eqnarray}\label{form4}
\Bigl.\frac{\partial_{z_1}^{l_j}\partial_{z_2}^{m_j}\psi_j^{q_j}}{l_j!\,m_j!}\Bigr|_{z_1={z_2}=0}=
\sum_{\begin{array}{c}\scriptstyle \lambda_1+\ldots+\lambda_{q_j}=l_j\\ \scriptstyle \mu_1+\ldots+\mu_{q_j}=m_j\end{array}}
(\lambda_1r_1+\mu_1r_2)^j\ldots(\lambda_{q_j}r_1+\mu_{q_j}r_2)^j\,c_{\lambda_1,\mu_1}\ldots c_{\lambda_{q_j},\mu_{q_j}}.
\end{eqnarray}

The summands in the right hand side of (\ref{form4}) do not contain the coefficient $c_{l,m}$. Indeed, if $(\lambda_i,\mu_i)=(l,m)$ in some summand,
this would necessary imply $(l_j,m_j)=(l,m)$ and, therefore, $(\lambda_s,\mu_s)=(0,0)$ for the other $(\lambda_s,\mu_s)$ in this summand. Thus, the formula
(\ref{clm}) can be written in the form
$$
L(lr_1+mr_2)\,c_{l,m}=p_{l,m}(\{\alpha_{p,q,Q}\},\{c_{\lambda,\mu}\}),
$$
where $p_{l,m}$ is the polynomial of the variables $\{\alpha_{p,q,Q}\}$, $\{c_{\lambda,\mu}\}$ (with $(\lambda,\mu)\in M$, $\lambda\leqslant l$, $\mu\leqslant m$, 
$(\lambda,\mu)\ne(l,m)$) determined by the formulae (\ref{form3}), (\ref{form4}).

Now we similarly use the equality (\ref{majorant}) to obtain recursive expressions for the coefficients $C_{l,m}$. Denote by $\Phi$ the power series from
the right hand side of this equality, 
$$
\Phi=\sum_{p,q,Q}|\alpha_{p,q,Q}|\,z_1^p{z_2}^qW^{q_0}W^{q_1}\ldots W^{q_n}\in{\mathbb C}\{z_1,{z_2}\},
$$
for $W=\sum\limits_{(l,m)\in{\mathbb Z}_+^2\setminus\{0\}}|lr_1+mr_2|^n\,C_{l,m}z_1^lz_2^m$. Then (\ref{majorant}) implies
\begin{eqnarray}\label{Clm}
\sigma\,|lr_1+mr_2|^n\,C_{l,m}=\Bigl.\frac{\partial_{z_1}^l\partial_{z_2}^m\Phi}{l!\,m!}\Bigr|_{z_1={z_2}=0}.
\end{eqnarray}
Keeping in mind the analogy and difference between the series $\phi\in{\mathbb C}[[z_1,{z_2}]]$ and $\Phi\in{\mathbb C}\{z_1,{z_2}\}$, we obtain
\begin{eqnarray}\label{form5}
\Bigl.\frac{\partial_{z_1}^l\partial_{z_2}^m\Phi}{l!\,m!}\Bigr|_{z_1={z_2}=0}=\sum_{p,q,Q}|\alpha_{p,q,Q}|
\sum_{\begin{array}{c}\scriptstyle l_0+l_1+\ldots+l_n=l-p\\ \scriptstyle m_0+m_1+\ldots+m_n=m-q\end{array}}
\Bigl.\frac{\partial_{z_1}^{l_0}\partial_{z_2}^{m_0}W^{q_0}}{l_0!\,m_0!}\;\ldots\;
\frac{\partial_{z_1}^{l_n}\partial_{z_2}^{m_n}W^{q_n}}{l_n!\,m_n!}\Bigr|_{z_1={z_2}=0},
\end{eqnarray}
where
\begin{eqnarray}\label{form6}
\Bigl.\frac{\partial_{z_1}^{l_j}\partial_{z_2}^{m_j}W^{q_j}}{l_j!\,m_j!}\Bigr|_{z_1={z_2}=0}=
\sum_{\begin{array}{c}\scriptstyle \lambda_1+\ldots+\lambda_{q_j}=l_j\\ \scriptstyle \mu_1+\ldots+\mu_{q_j}=m_j\end{array}}
|\lambda_1r_1+\mu_1r_2|^n\ldots|\lambda_{q_j}r_1+\mu_{q_j}r_2|^n\,C_{\lambda_1,\mu_1}\ldots C_{\lambda_{q_j},\mu_{q_j}}.
\end{eqnarray}
Thus, the formula (\ref{Clm}) can be written in the form
\begin{eqnarray}\label{Plm}
\sigma\,|lr_1+mr_2|^n\,C_{l,m}=P_{l,m}(\{|\alpha_{p,q,Q}|\},\{C_{\lambda,\mu}\}),
\end{eqnarray}
where $P_{l,m}$ is the polynomial of the variables $\{|\alpha_{p,q,Q}|\}$, $\{C_{\lambda,\mu}\}$ ($\lambda\leqslant l$, $\mu\leqslant m$, $(\lambda,\mu)\ne(l,m)$)
with the real positive coefficients determined by the formulae (\ref{form5}), (\ref{form6}). Since for $(l,m)$ equal to $(1,0)$ and $(0,1)$ we have
$$
\sigma\,|r_1|^n\,C_{1,0}=\partial_{z_1}\Phi(0,0)=|\alpha_{1,0,{\bf 0}}|, \qquad \sigma\,|r_2|^n\,C_{0,1}=\partial_{z_2}\Phi(0,0)=|\alpha_{0,1,{\bf 0}}|,
$$
all the coefficients $C_{l,m}$ are real nonnegative numbers.

Finally we come to a conclusive part of the proof, the estimates
$$
|c_{l,m}|\leqslant C_{l,m}\quad\forall (l,m)\in{\mathbb Z}_+^2\setminus\{0\}.
$$
We prove them by the induction with respect to the sum $l+m$ of the indices.

For $l+m=1$ according to (\ref{clm}) we have
$$
L(r_1)\,c_{1,0}=\partial_{z_1}\phi(0,0)=\alpha_{1,0,{\bf 0}}, \qquad L(r_2)\,c_{0,1}=\partial_{z_2}\phi(0,0)=\alpha_{0,1,{\bf 0}},
$$
hence\footnote{Note that if $(1,0)\in M$, then ${|L(r_1)|}/{|r_1|^n}\geqslant\sigma$, and in the opposite case $c_{1,0}=0$. The similar is true for
the index $(0,1)$.}
$$
|c_{1,0}|=\frac{|\alpha_{1,0,{\bf 0}}|}{|L(r_1)|}=\frac{\sigma\,|r_1|^n\,C_{1,0}}{|L(r_1)|}\leqslant C_{1,0}, \qquad
|c_{0,1}|=\frac{|\alpha_{0,1,{\bf 0}}|}{|L(r_2)|}=\frac{\sigma\,|r_2|^n\,C_{0,1}}{|L(r_2)|}\leqslant C_{0,1}.
$$
Further, by the construction of the polynomials $p_{l,m}$ and $P_{l,m}$, for any $(l,m)\in M$ we have
$$
\bigl|p_{l,m}(\{\alpha_{p,q,Q}\},\{c_{\lambda,\mu}\})\bigr|\leqslant P_{l,m}(\{|\alpha_{p,q,Q}|\},\{|c_{\lambda,\mu}|\})
$$
(here we use the estimate $|\lambda r_1+\mu r_2|\geqslant1$ for all $(\lambda,\,\mu)\in M$), and the inductive assumption (the second inequality below) implies
$$
|L(lr_1+mr_2)|\,|c_{l,m}|=\bigl|p_{l,m}(\{\alpha_{p,q,Q}\},\{c_{\lambda,\mu}\})\bigr|\leqslant P_{l,m}(\{|\alpha_{p,q,Q}|\},\{|c_{\lambda,\mu}|\})
\leqslant $$
\begin{equation}\label{baikal} 
\leqslant P_{l,m}(\{|\alpha_{p,q,Q}|\},\{C_{\lambda,\mu}\})\leqslant\sigma\,|lr_1+mr_2|^n\,C_{l,m},\footnote{In the left hand side of the last inequality in 
(\ref{baikal}) the indices $(\lambda,\mu)$ of the variables $C_{\lambda,\mu}$ belong to the set $M$, while in the equality (\ref{Plm}) they belong to 
${\mathbb Z}_+^2\setminus\{0\}$. Therefore we write in (\ref{baikal}) the inequality instead of the equality.}
\end{equation}
whence the required estimates follow:
$$
|c_{l,m}|\leqslant\frac{\sigma\,|lr_1+mr_2|^n}{|L(lr_1+mr_2)|}\,C_{l,m}\leqslant C_{l,m}
$$
(for $(l,m)\in M$ we have $|L(lr_1+mr_2)|/|lr_1+mr_2|^n\geqslant\sigma$, and $c_{l,m}=0$ for $(l,m)\not\in M$).

Now it remains to note that for any sector $S$ with the vertex at $0\in\mathbb C$ and of the opening less than $2\pi$, the terms of the series (\ref{series}) 
are regarded as holomorphic single-valued functions in $S$, and to pass from the convergence of $\psi_0=\sum_{(l,m)\in{\mathbb Z}_+^2\setminus\{0\}}c_{l,m}z_1^lz_2^m$ to the convergence of $\varphi=\sum_{k=0}^{\infty}c_kx^{s_k}$. Let the power series $\psi_0$ converge in a neighbourhood of a closed polydisk 
$\{|z_1|\leqslant R,\,|z_2|\leqslant R\}$. Then there is a positive constant $C$ such that
$$
|c_{l,m}|\leqslant C/R^{l+m}
$$
for all $(l,m)\in{\mathbb Z}_+^2\setminus\{0\}$. If $x\in S$ is small enough for the inequalities
$$
|x^{r_1}|=|x|^{{\rm Re}\,r_1}e^{-{\rm Im}\,r_1\cdot\arg x}\leqslant R_1<R, \qquad |x^{r_2}|=|x|^{{\rm Re}\,r_2}e^{-{\rm Im}\,r_2\cdot\arg x}\leqslant R_1<R
$$
to be held (recall that ${\rm Re}\,r_{1,2}>0$), then
$$
|c_kx^{s_k-s_{\mu}}|=|c_{l,m}|\cdot|x^{r_1}|^l\cdot|x^{r_2}|^m\leqslant C\left(\frac{R_1}R\right)^l\left(\frac{R_1}R\right)^m=C\,h^{l+m},\quad l=l(k),\,m=m(k).
$$
As $0<h<1$, the series $\sum_{(l,m)\in{\mathbb Z}_+^2}h^{l+m}$ converges, hence the series $\varphi=\sum_{k=0}^{\infty}c_kx^{s_k}$ converges uniformly in $S$ for
sufficiently small $|x|$. 
\medskip

{\bf Acknowledgements.} We are very thankful to Anton A.\,Vladimirov who has kindly proved Lemma 3 for us.

Renat Gontsov

Institute for Information Transmission Problems, Moscow, Russia

rgontsov@inbox.ru
\\
\\

Irina Goryuchkina

Institute of Applied Mathematics, Moscow, Russia

igoryuchkina@gmail.com


\begin{thebibliography}{99}

\bibitem{GS}
D.\,Yu.\,Grigor'ev, M.\,F.\,Singer, {\it Solving ordinary differential equations in terms of series with real exponents},
Trans. Amer. Math. Soc., V.~327(1) (1991), 329--351.

\bibitem{Br}
A.\,D.\,Bruno, {\it Asymptotic behaviour and expansions of solutions of an ordinary differential equation}, Russian Math. Surv.,
V.~59(3) (2004), 429--480.

\bibitem{Mal}
B.\,Malgrange, {\it Sur le th\'eor\`eme de Maillet}, Asympt. Anal., V.~2 (1989), 1--4.

\bibitem{GG}
R.\,R.\,Gontsov, I.\,V.\,Goryuchkina, {\it An analytic proof of the Malgrange--Sibuya theorem on the convergence of formal solutions
of an ODE}, arXiv:1311.6416(math.CA), 9 pp.

\bibitem{BG}
A.\,D.\,Bruno, I.\,V.\,Goryuchkina, {\it Asymptotic expansions of the solutions of the sixth Painlev\'e equation},
Trans. Moscow Math. Soc. (2010), 1--104.

\bibitem{Ca}
J.\,Cano, {\it On the series defined by differential equations, with an extension of the Puiseux polygon construction to these
equations}, Analysis, V.~13 (1993), 103--119.

\end{thebibliography}
\end{document}